\documentclass[11pt]{article}

\usepackage{amssymb}
\usepackage{amsthm}
\usepackage{amsmath,amsfonts,amssymb}

\allowdisplaybreaks
\usepackage[mathscr]{eucal}
\usepackage{exscale}
\usepackage[square,sort,comma,numbers]{natbib}
\setcitestyle{authoryear,round}

\usepackage{bm}
\usepackage{eqlist}
\usepackage[dvipsnames]{color}
\usepackage[Lenny]{fncychap}
\usepackage{multirow}
\usepackage{graphicx}
\usepackage{subfigure}
\usepackage{algpseudocode}
\usepackage{algorithm}
\usepackage{caption}
\usepackage{hyperref}
\usepackage{framed}
\usepackage{color}
\usepackage{booktabs}
\usepackage{makecell}
\usepackage{threeparttable}
\usepackage{setspace}
\usepackage{framed}
\usepackage{algpseudocode}
\usepackage{indentfirst} 
\usepackage{longtable}
\usepackage{tabu}
\usepackage{tikz}
\usepackage{chngpage}

\usepackage[colorinlistoftodos]{todonotes}

\usepackage{ulem}
\usepackage{soul}

\usepackage[title]{appendix}
\usepackage[utf8]{inputenc} % this is needed for umlauts
\usepackage[english]{babel} % this is needed for umlauts
\usepackage[T1]{fontenc}    % this is needed for correct output of umlauts in pdf
\usepackage{amssymb,amsmath,amsfonts} % nice math rendering
\usepackage{braket} % needed for \Set
\usepackage{caption}
\usepackage{algorithm}

\hypersetup{
	colorlinks=true,
	linkcolor=blue,
	citecolor=blue,
	citebordercolor={1 1 0},
}

\usepackage{setspace}

\font\bigbf=cmbx10 scaled \magstep3

\title{\bigbf  Enhancing the sensing power of bike-sharing system for urban environment}

\author{Wen Ji \textsuperscript{a,c,1},
	Ke Han \textsuperscript{b,1},
	Qi Hao \textsuperscript{c},
	Qian Ge \textsuperscript{a},
	Ying Long\textsuperscript{c,d} \thanks{Corresponding author} \\
	\small \textsuperscript{a} \textit{School of Transportation and Logistics,} \\
	 \textit{Southwest Jiaotong University, Chengdu, Sichuan 611756, China} \\
	\small \textsuperscript{b} \textit{School of Economics and Management,}\\
	 \textit{Southwest Jiaotong University, Chengdu, Sichuan 610031, China} \\
	\small \textsuperscript{c} \textit{School of Architecture,} \\
	\textit{ Tsinghua University, Beijing, 100084, China} \\
	\small \textsuperscript{d} \textit{Hang Lung Center for Real Estate, Key Laboratory of Ecological Planning $\&$ Green Building,} \\ \textit{Ministry of Education, Tsinghua University, Beijing, 100084, China}
}
\date{}

\begin{document}

	\maketitle
	\renewcommand{\thefootnote}{1}
	\footnotetext[1]{These authors contributed equally to this work.}
	\renewcommand{\thefootnote}{}
	\footnotetext{E-mail addresses: wenji@my.swjtu.edu.cn (W. Ji), kehan@swjtu.edu.cn (K. Han), hao-q14@tsinghua.org.cn (Q. Hao), geqian@swjtu.edu.cn (Q. Ge), ylong@tsinghua.edu.cn (Y. Long)}
	
	\begin{abstract}
		The development of smart cities requires innovative sensing solutions for efficient and low-cost urban environment monitoring. Bike-sharing systems, with their wide coverage, flexible mobility, and dense urban distribution, present a promising platform for pervasive sensing. At a relative early stage, research on bike-based sensing focuses on the application of data collected via passive sensing, without consideration of the  optimization of  data collection through sensor deployment or vehicle scheduling. To address this gap, this study integrates a binomial probability model with a mixed-integer linear programming model to optimize sensor allocation across bike stands. Additionally, an active scheduling strategy guides user bike selection to enhance the efficacy of data collection. A case study in Manhattan validates the proposed strategy, showing that equipping sensors on just 1\% of the bikes covers approximately 70\% of road segments in a day, highlighting the significant potential of bike-sharing systems for urban sensing.
	\end{abstract}
	
	
	\section{Introduction}\label{secIntro}
	The rapid urbanization of cities has created a pressing need for innovative urban sensing solutions to monitor and understand complex urban dynamics efficiently \citep{DSXVF2019}. Among the various sensing paradigms, drive-by sensing (DS) has emerged as a particularly promising approach, leveraging the mobility of sensor-equipped vehicles to achieve extensive spatial-temporal coverage at relatively low costs \citep{MZY2014}. DS has been effectively applied in diverse scenarios, including air quality \citep{HSWHFABT2015, SM2019}, urban heat island phenomena \citep{FBRDRAP2021}, noise \citep{Alsina-Pages2017}, traffic conditions estimation \citep{AQSS2021}, urban greenery \citep{GMPDPR2024}, storefront vacancy \citep{LL2024} and abandoned buildings \citep{LMZLL2023}. 
	
	Traditionally, DS has relied on various types of vehicles, such as taxis \citep{C2020, XCPJZN}, buses \citep{JHL2023b, AWSDR2024}, and dedicated vehicles \citep{Messier2018}. However, these vehicles face limitations in coverage, route constraints, and high deployment costs. Building on this foundation, bike-sharing systems have emerged as a promising platform for urban sensing due to their wide coverage, flexible mobility, and dense distribution across cities \citep{SLZSW2022}. While still in its early stages, research has begun to leverage bike-sharing systems for data collection. For instance,  \cite{WXLLWZLG2020} and \cite{VBP2020} both use bikes equipped with sensors to collect air quality data, demonstrating the feasibility of bike-based sensing paradigm. However, these studies all focus on applications of data collected via passive sensing for urban research, giving no consideration to the potential optimization of  data collection through sensor deployment or vehicle scheduling, leaving a significant gap in fully realizing their powers as urban sensing hosts. 
	
	In this context, our study systematically enhances the sensing power of bike-sharing systems by addressing both strategic and operational aspects. Strategically, we optimize sensor allocation across bike stands to achieve maximum spatial coverage, while operationally, we introduce an active scheduling approach to guide user bike selection, thereby improving the data collection efficiency of sensor-equipped bikes. The main contributions of this paper are as follows:
	
	\begin{itemize}
		\item We introduce a binomial probability model to capture the coverage of road segments by sharing bikes. Then, we formulate the sensor-to-stand allocation problem as a Mixed-Integer Linear Programming (MILP) model to strategically maximize spatial coverage.

		\item We propose an active scheduling approach that leverages user bike selection guidance to enhance the sensing rewards of sensor-equipped bikes.
		
		\item A case study conducted in Manhattan validates the effectiveness of the proposed sensor optimization strategy and quantitatively assesses the spatial-temporal sensing power of the bike-sharing system. 
	\end{itemize}
	
	The remainder of this paper is organized as follows: Section \ref{sec_RW} provides a review of relevant literature. Section \ref{sec_method} outlines the research framework, methodology, and data employed in this study. The results and corresponding discussions are presented in Section \ref{sec_result}, and finally, the conclusions are drawn in Section \ref{sec_conclusion}.
	
	\section{Related work} \label{sec_RW}
	
	\subsection{Drive-by sensing}
	DS has emerged as a promising urban sensing paradigm, gaining significant attention for its ability to provide extensive spatial-temporal coverage at relatively low costs by leveraging the mobility of sensor-equipped vehicles \citep{BHK2021, DSXVF2019, JHL2023a}. However, the spatial-temporal sampling quality of DS is largely influenced by the mobility patterns of the host vehicles \citep{ADRMdR2018}. For instance, taxis, driven by profit-oriented operations, tend to cluster in high-demand areas, leading to spatial sampling biases \citep{OASSR2019, C2020}. Buses, restricted to fixed routes, are incapable of collecting data beyond their network coverage \citep{JHL2023b, AWSDR2024, HLYK2024}. Meanwhile, dedicated sensing vehicles, due to their high deployment and operational costs, are not suitable for large-scale deployment and are typically used only for targeted investigations at specific locations \citep{Messier2018, GM2020}.
	
	In contrast, bike-sharing systems have become increasingly prevalent in urban areas, with extensive coverage and dense distribution of bike stands, making them a promising platform for fine-grained urban sensing. Despite their potential, the sensing power of bike-sharing systems has not yet been explored in the existing literature.
	
	\subsection{Sensor deployment problems in DS}
	
	Research on sensor deployment for DS mainly focuses on taxis and buses. These studies typically assume that vehicle trajectories are known in advance and aim to select an optimal subset of vehicles for sensor installation to maximize sensing rewards. The problem is often formulated as subset selection problems, for which various heuristic algorithms are developed to provide effective solutions. For instance, \cite{ZML2014} design a greedy heuristic algorithm to select the minimum number of taxis required to meet specified monitoring requirements. \cite{ZYL2018} develop a multi-objective genetic algorithm to identify an optimal subset of taxis that jointly optimizes the overall coverage of all vehicles and the reliability of individual sensing tasks. \cite{Kaivonen2020} employs a route coverage image analysis algorithm to determine the best combination of bus routes for maximizing spatial coverage. \cite{Gao2016} design a greedy heuristic algorithm to select a subset of buses that maximizes the sensing reward. \cite{Agarwal2020} introduce an approximate algorithm for selecting a subset of vehicles from a mixed fleet of buses and taxis to maximize sensing reward.
	
	Different from the studies mentioned above, bike trajectories are uncertain, making it challenging to apply similar methods directly. This study first employs a binomial probability model to estimate the relationship between the number of bikes at stands and the coverage of road segments. The sensor-stand allocation problem is then formulated as an MILP model, which can be directly solved by off-the-shelf solvers.
	
	\subsection{Active scheduling of sensor-equipped vehicles}
	
	The potential of actively scheduling sensor-equipped vehicles to enhance sensing power attracts significant research interest. For taxis, strategies such as optimizing vehicle-passenger matching \citep{MH2023, C2020} and incentivizing routing for both occupied \citep{ADFS2021} and unoccupied vehicles \citep{GQ2024} are employed to improve sensing rewards. In the case of buses, some studies focus on the joint optimization of timetable coordination and vehicle scheduling to maximize sensing rewards \citep{JHL2023b, DH2023}. For dedicated vehicles, specific route planning algorithms are typically developed to achieve maximum sensing efficiency \citep{JHG2023}.
	
	This study investigates the potential to enhance the sensing power of bike-sharing systems by guiding users' bike selection behavior.
	
	\subsection{Bike-sharing system}
	Bike-sharing systems have garnered significant attention as an innovative urban transportation mode in recent years. Current research primarily concentrates on demand prediction \citep{LHP2018, XJL2018}, rebalance scheduling \citep{CMLZH2024, HYCCC2024}, operational policies \citep{KJS2020, YXLZLSCYY2022}, travel behavior \citep{YHTC2019, KFDDSR2022}, and environmental impact \citep{LZB2018, ZSLXFZHSL2019}.
	
	However, beyond their role as transportation tools, bike-sharing systems can be envisioned as mobile sensors capable of collecting valuable data on urban environments and dynamics. Exploring this perspective could provide new insights into the application scenarios of bike-sharing systems, thereby presenting significant academic value and practical relevance.
	
	\section{Material and Methods} \label{sec_method}
	
	\subsection{Overview}
\begin{figure}[h!]
		\centering
		\includegraphics[width=1.0\textwidth]{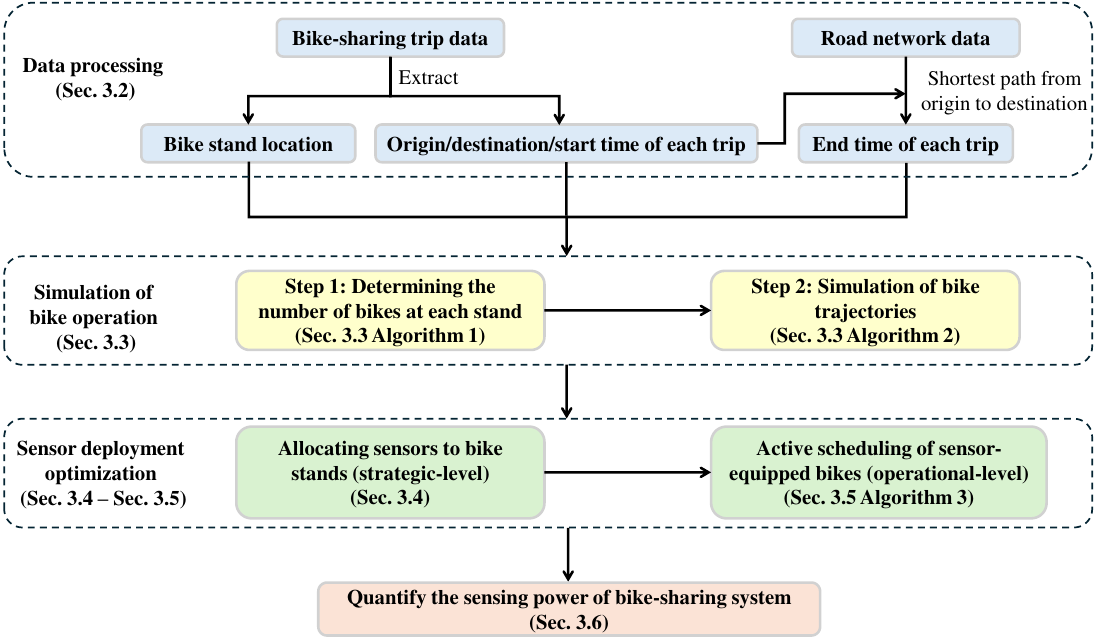}
		\caption{Methodological framework.}
		\label{fig_overview}
	\end{figure}

	This study employs a data-driven methodology framework to optimize sensor deployment and the active scheduling of sensor-equipped bikes, aiming to maximize data collection efficiency. The methodological framework is shown in Figure \ref{fig_overview}. Table \ref{tab_notations} lists key parameters and variables used in this paper.
	
	\begin{itemize}
		\item[(1)] \textbf{Data processing}: The input includes bike-sharing trip data and road network information, with the process focusing on cleaning the trip data and extracting bike stand locations. Further details are provided in Section \ref{sec_dp}.
		\item[(2)] \textbf{Simulation of bike operation}: This part introduces a simulation method based on bike-sharing trip data and stand information to estimate the total number of bikes and generate individual bike trajectories. Further details can be found in Section \ref{sec_simulation}.
		\item[(3)] \textbf{Allocating sensors to bike stands}: This part addresses the allocation of sensors to bike stands to maximize spatial coverage by formulating an MILP model. Further details are provided in Section \ref{sec_MILP}.
		\item[(4)] \textbf{Active scheduling of sensor-equipped bikes}: Building on the results of sensor allocation, this part proposes an active scheduling strategy for sensor-equipped bikes to enhance the sensing reward further. Further details are provided in Section \ref{sec_active}.
		\item[(5)] \textbf{Quantify the sensing power of bike-sharing system}: Utilizing the simulation optimization method described above, this part introduces a metric to quantify the sensing power of the bike-sharing system under various monitoring requirements. Further details are provided in Section \ref{sec_quantify}.
	\end{itemize}

	\setlength\LTleft{0pt}
	\setlength\LTright{0pt}
	\begin{longtable}{@{\extracolsep{\fill}}rl}
		\caption{Notations and symbols}
		\label{tab_notations} 
		\\
		\hline
		\multicolumn{2}{l}{Sets}     
		\\
		\hline
		$\mathcal{R}$  & Set of bike trips;
		\\
		$\mathcal{N}$  & Set of nodes in road network;
		\\
		$\mathcal{A}$  & Set of road segments in road network;
		\\
		$\mathcal{T}$  & Set of discrete time intervals;
		\\
		$\mathcal{S}$  & Set of bike stands, $\mathcal{S} \subseteq \mathcal{N}$;
		\\
		$\mathcal{K}_{s}$  & Set of bikes at each stand $s$;
		\\
		\hline
		\multicolumn{2}{l}{Parameters and constants}    
		\\  
		\hline
		$b_{s}$ & Initial number of bikes deployed at stand $s$, $s \in \mathcal{S}$;
		\\
		$t_0$  & Simulation start time;
		\\
		$T$  & Simulation end time;
		\\
		$N_{s,e}$ & Number of times the road section $e$ has been covered by sensor-equipped bikes
		\\
		& at stand $s$ during simulation time;
		\\
		$N_s$  & Total number of sensors;
		\\
		$N_{\mathcal{T}}$ & Total number of time intervals;
		\\
		$p_{s,e}$ & Probability that a bike from stand $s$ intersects road segment $e$ during simulation
		\\
		& time;
		\\
		$\beta$ & Probability that a user accepts the guidance to choose a sensor-equipped bike;
		\\
		$l_e$  & Length of road segment $e$;
		\\
		$N_{e,t}$  & Number of sensor-equipped bikes covering road segment $e$ during time interval $t$;
		\\
		\hline
		\multicolumn{2}{l}{Auxiliary variables}                
		\\
		\hline
		$y_{e}$  & Binary variable that equals 1 if the expected coverage of road segment $e$ exceed 
		\\
		& $K$ times;
		\\
		$N_{e}$  & Number of times the road segment $e$ is covered;
		\\
		\hline
		\multicolumn{2}{l}{Decision variables}                 
		\\
		\hline
		$n_{s}$  &  Number of sensors allocated to bike stand $s$.
		\\
		\hline
	\end{longtable}
	
	\subsection{Data processing} \label{sec_dp}
\begin{figure}[h!]
		\centering
		\includegraphics[width=0.9\textwidth]{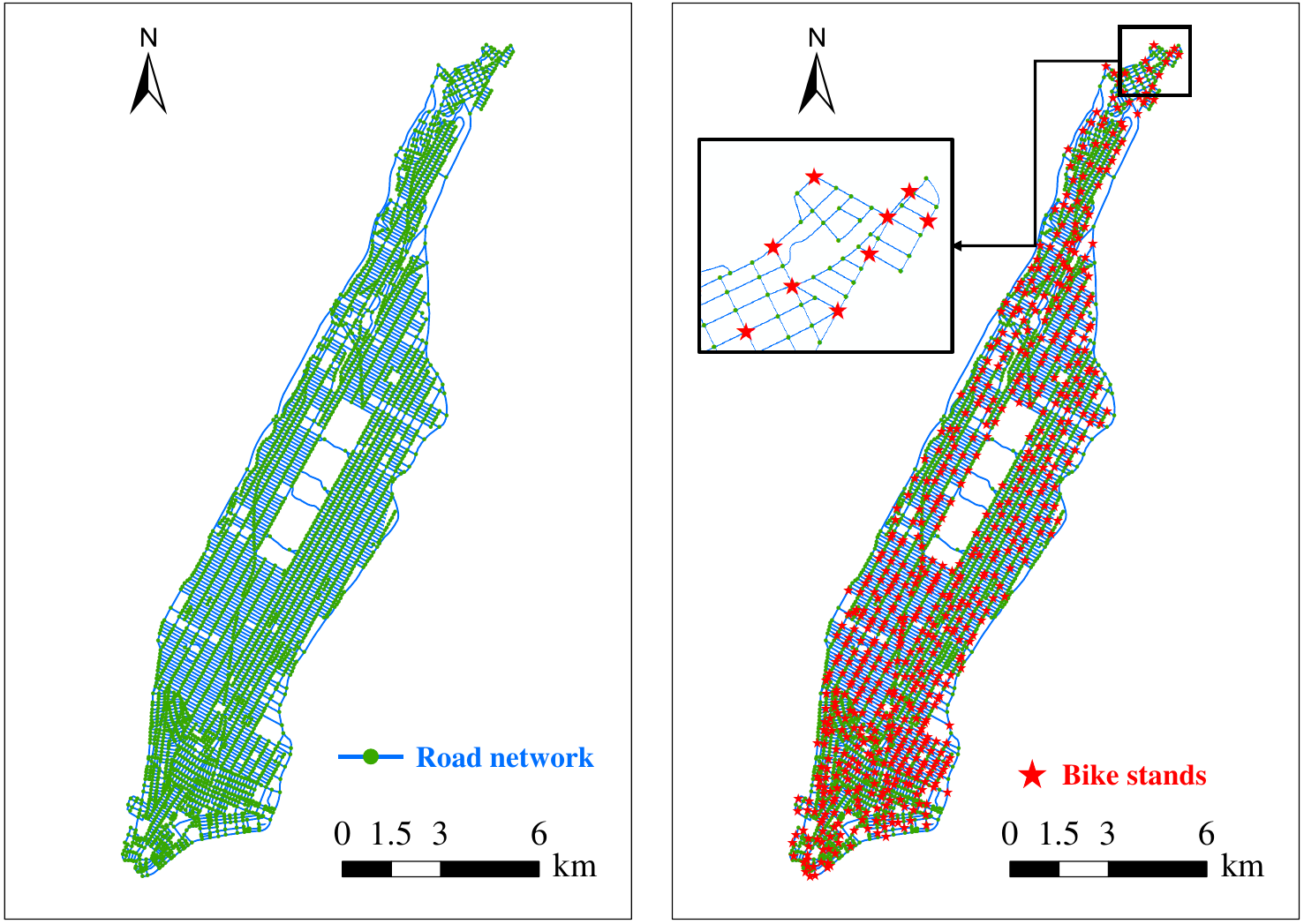}
		\caption{Left: Road network of Manhattan, New York City, which consists of 3,330 nodes and 6,055 road segments. Right: Locations of 646 bike-sharing stands in Manhattan.}
		\label{fig_network}
	\end{figure}

	The bike-sharing trip data for Manhattan are obtained from the official website of Citi Bike operators\renewcommand{\thefootnote}{2} \footnote{NYC (Citi): \url{https://www.citibikenyc.com/system-data}}. We use trip data from March 1, 2024, to extract the origin and destination stand locations, as well as the start time for each trip. The Manhattan road network data are retrieved from the OpenStreetMap website \renewcommand{\thefootnote}{3}\footnote{OpenStreetMap: \url{https://www.openstreetmap.org}},  as shown in Figure \ref{fig_network} (Left). The bike-sharing trip data are then preprocessed according to the following steps:
	
	\begin{itemize}
		\item \textbf{Step 1:} The bike stands are mapped to the nearest nodes on the road network based on their coordinates. This process identifies 646 unique bike stands within the Manhattan road network, as shown in Figure \ref{fig_network} (Right).
		\item \textbf{Step 2:} We use the road network adjacency matrix to calculate the shortest path distance between the origin and destination of each trip. Assuming a constant bike speed of 13 km/h \citep{KFDDSR2022}, we estimate the trip duration and end time \renewcommand{\thefootnote}{4}\footnote{Since the original data for sharing bikes does not include detailed travel routes, using the reported end time of each trip does not allow for accurate trajectory reconstruction.}.
		\item \textbf{Step 3:} We retain only trips with a total distance between 0.5 km and 5 km, a common data cleaning step employed in many studies \citep{ZLL2019, YHLLCKL2024}.
		\item \textbf{Step 4:} Trips are further filtered to include only those with start times between 6 am and 10 pm, as the demand for data collection during nighttime hours is minimal.
	\end{itemize}

	\subsection{Simulation of bike operations based on trip data} \label{sec_simulation}
	This section introduces a method for simulating bike trajectories based on trip data. It primarily involves determining the initial number of bikes at each stand and simulating the trajectory of each bike.
	
	When determining the initial number of bikes at each stand, it is essential to ensure that all trips can be served. This is achieved by analyzing the traffic flow at each stand throughout the simulation time period. At each time step, the number of bikes for each stand $s$ is updated and recorded based on the trip data. The detailed process is outlined in Algorithm \ref{alg_bike_num}.
	
	 \begin{algorithm}[h!]
		\caption{Determining the initial number of bikes at each stand}
		\begin{tabbing}
			\hspace{0.01 in}\=  \hspace{0.9 in}\= \kill % set up two tab positions
			\>{\bf Input}  \> Set of bike trips $\mathcal{R}$, set of bike stands $\mathcal{S}$, simulation start time $t_0$, 
			\\
			\>  \> simulation end time $T$;
			\\
			\> {\bf Initialize} \> Initial number of bikes deployed at stand at start time $t_{0}$: $b_s^{0}=0, \forall s \in \mathcal{S}$, 
			\\
			\> \> list of the number of bikes at each stand $s$ at each time step: 
			\\
			\> \> $\phi_s = [b_s], \forall s \in \mathcal{S}$, $t = t_0$, $b_s^t = b_s^{0}$;
			\\
			\> {\bf 1:} \> For $\forall s \in \mathcal{S}$:
			\\
			\> {\bf 2:} \> \hspace{0.2 in} While $t \leq T$:
			\\
			\> {\bf 3:} \> \hspace{0.4 in}  Count the number of arriving trips at stand $s$ at time $t$: $n_{arr, t}^s$
			\\
			\> {\bf 4:} \> \hspace{0.4 in} Count the number of departing trips from stand $s$ at time $t$: $n_{dep, t}^s$
			\\
			\> {\bf 5:} \> \hspace{0.4 in} Update the number of bikes at stand $s$ at time $t$: 
			\\
			\> \> \hspace{0.4 in} $b_s^t \leftarrow b_s^t + n_{arr, t}^s - n_{dep, t}^s$
			\\
			\> {\bf 6:} \> \hspace{0.4 in} Record the number of bikes at stand $s$ at time $t$: append $b_s^t$ to $\phi_s$
			\\
			\> {\bf 7:} \> \hspace{0.4 in} Update time: $t = t + 1$;
			\\
			\> {\bf 8:} \> \hspace{0.2 in} End While
			\\
			\> {\bf 9:} \> \hspace{0.2 in} Determine the initial number of bikes to be deployed at stand $s$:
			\\
			\> \> \hspace{0.2 in} $b_s = b_s^{0} - \min(\phi_s)$;
			\\
			\> {\bf 10:} \> End For
			\\
			\>{\bf Output}      \> Initial number of bikes at each stand $b_{s}, \forall s \in \mathcal{S}$.
		\end{tabbing}
		\label{alg_bike_num}
	\end{algorithm}
	
	After determining the initial number of bikes at each stand, we denote the set of bikes at each stand $s$ as $\mathcal{K}_{s}$. The next step involves tracking the trajectory of each bike, specifically the sequence of trips it serves, using real-time trip data. The core idea of the simulation process is to select a random bike from the available bikes at the trip's starting stand for each trip, thereby reflecting real-world conditions. The detailed simulation process is outlined in Algorithm \ref{alg_bike_trajectory}.
	
	\begin{algorithm}[h!]
		\caption{Simulation of bike trajectories}
		\begin{tabbing}
			\hspace{0.01 in}\=  \hspace{0.9 in}\= \kill % set up two tab positions
			\>{\bf Input}  \> Set of bike trips $\mathcal{R}$, set of bike stands $\mathcal{S}$, simulation start time $t_0$, simulation 
			\\
			\>  \> end time $T$, initial number of bikes at each stand $b_{s}, \forall s \in \mathcal{S}$, set of bikes at 
			\\
			\> \> each stand $\mathcal{K}_s, \forall s \in \mathcal{S}$;
			\\
			\> {\bf Initialize} \> $t = t_0$, the status of each bike $k \in \mathcal{K}_s, \forall s \in \mathcal{S}$ is set to idle;
			\\
			\> {\bf 1:} \>  While $t \leq T$:
			\\
			\> {\bf 2:}\> \hspace{0.2 in} Check if any busy bikes need to end the trip at time $t$. If so, update their
			\\
			\> \>  \hspace{0.2 in}   status to idle;
			\\
			\> {\bf 3:}\>  \hspace{0.2 in} Find all trips starting at time $t$, denoted as $\mathcal{R}^{'}$;
			\\
			\> {\bf 4:}\>   \hspace{0.2 in} For $r \in \mathcal{R}^{'}$:
			\\
			\> {\bf 5:}\>   \hspace{0.4 in} Randomly select an idle bike from the start stand of trip $r$;
			\\
			\> {\bf 6:}\>  \hspace{0.4 in} Record the trajectory of the selected bike and update its status to busy;
			\\
			\> {\bf 7:}\> \hspace{0.2 in} End For
			\\
			\> {\bf 8:}\> \hspace{0.2 in} Update time: $t = t + 1$;
			\\
			\> {\bf 9:} \> End While
			\\
			\>{\bf Output}      \> The trajectory of each bike $k \in \mathcal{K}_s, \forall s \in \mathcal{S}$.
		\end{tabbing}
		\label{alg_bike_trajectory}
	\end{algorithm}
	
	\subsection{Allocating sensors to bike stands} \label{sec_MILP}
	
	In this section, we first explore the relationship between the number of bikes at stands and the coverage of road segments. Next, we address the sensor allocation problem by developing an MILP model. The objective is to determine the optimal distribution of sensors across various bike stands to maximize spatial coverage.
	
	\subsection*{Relationship between the number of bikes at stands and the coverage of road segments}
	
	\textbf{Hypothesis:} Inspired by \cite{OASSR2019} and \cite{HJNLL2024}, we hypothesize that the coverage of a road segment by bikes originating from a single stand can be modeled using a binomial distribution. Specifically, assuming that $n_s$ bikes, all from the same stand $s$, operate independently and follow a similar spatial distribution across the network, the coverage of road section $e \in \mathcal{A}$ by bikes at stand $s \in \mathcal{S}$ within the simulation time period $T-t_{0}$ is binomially distributed as $c_{s,e} \sim B(n_s, p_{s,e})$, where $p_{s,e}$ represents the probability that a bike from stand $s$ visits road segment $e$ during time period $T-t_{0}$. Consequently, during the simulation period, the average number of times bikes at stand $s$ cover road segment $e$ can be expressed as the expected value of a binomial distribution:
	
	\begin{equation}\label{eq_B1}
		N_{s,e} \approx p_{s,e}n_s \qquad \forall s\in \mathcal{S},\,e\in\mathcal{A}
	\end{equation}

	\textbf{Verification:} To validate the approximation in \eqref{eq_B1}, we employed the bike trajectory simulation method described in Algorithm \ref{alg_bike_trajectory}, using real-world bike-sharing trip data from Manhattan. Specifically, we analyzed 16 hours of data (6 am to 10 pm) from March 1, 2024. Given the inherent randomness of the simulation process outlined in Algorithm \ref{alg_bike_trajectory}, we conducted 20 independent simulation runs and averaged the results to obtain reliable estimates. For a fixed bike stand $s$ and a fixed number of bikes $n_s$, let $N_{s,e}(\tau)$ be the number of times the road section $e$ has been covered during 16 hours in the $\tau$-th simulation $(\tau=1, ..., 20)$, and define:
	
	\begin{equation}\label{eq_B2}
	\bar N_{s,e}={1\over 20}\sum_{\tau=1}^{20}N_{s,e}(\tau)\qquad \forall s\in \mathcal{S}, e \in \mathcal{A}
	\end{equation}
	
\begin{figure}[h!]
		\centering
		\includegraphics[width=1.0\textwidth]{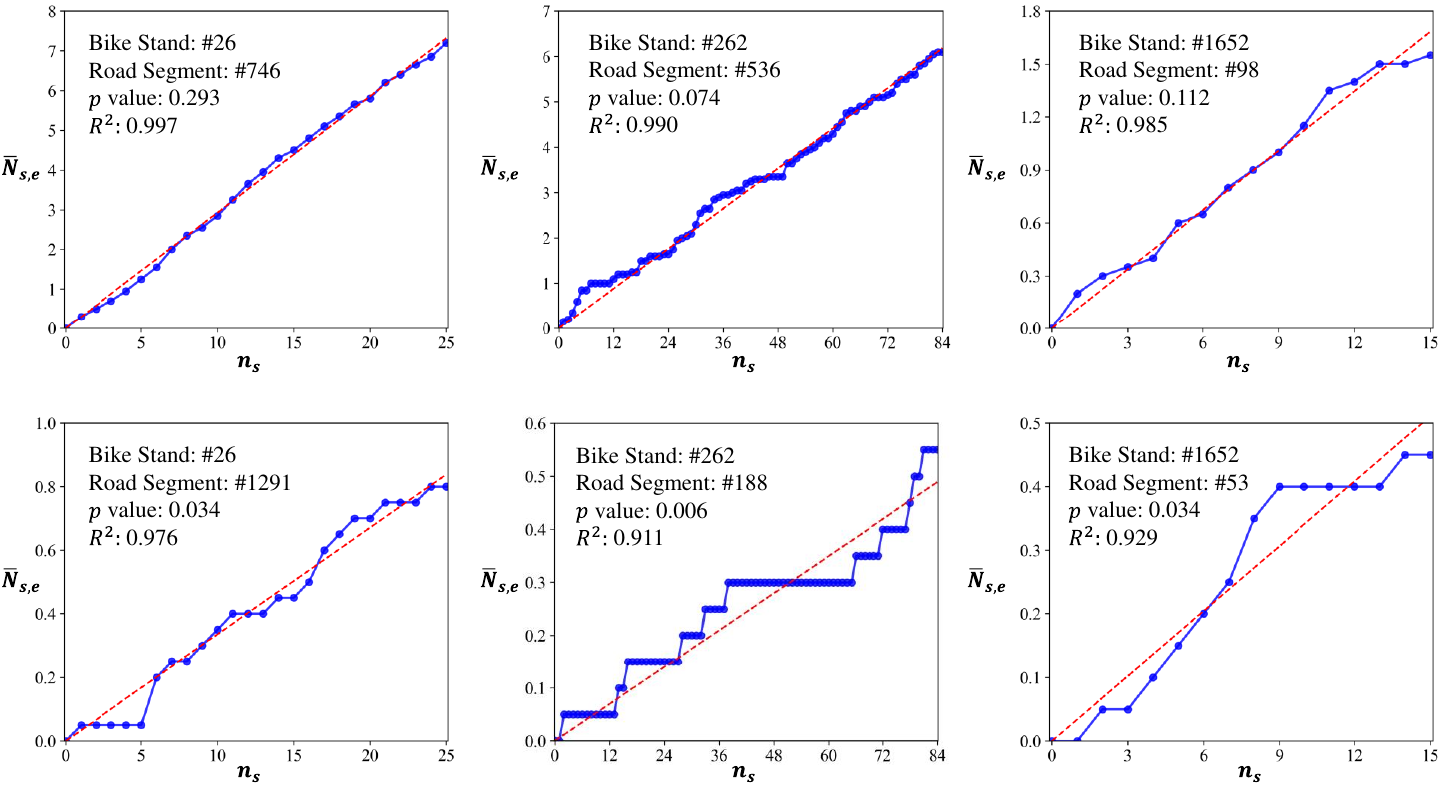}
		\caption{Empirical relationship between $n^s$ and average coverage $\bar N_{s,e}$.}
		\label{fig_bi}
	\end{figure}

	The empirical relationships between $n_s$ and $\bar N_{s,e}$ are illustrated in Figure \ref{fig_bi}. We see that, for a fixed bike stand $s$ and a fixed road segment $e$, the relationships are approximately linear, which agrees with \eqref{eq_B1}. Furthermore, it can be observed that the fitting effect improves with the increase of $n_s$ and $\bar N_{s,e}$.

\begin{figure}[h!]
		\centering
		\includegraphics[width=1.0\textwidth]{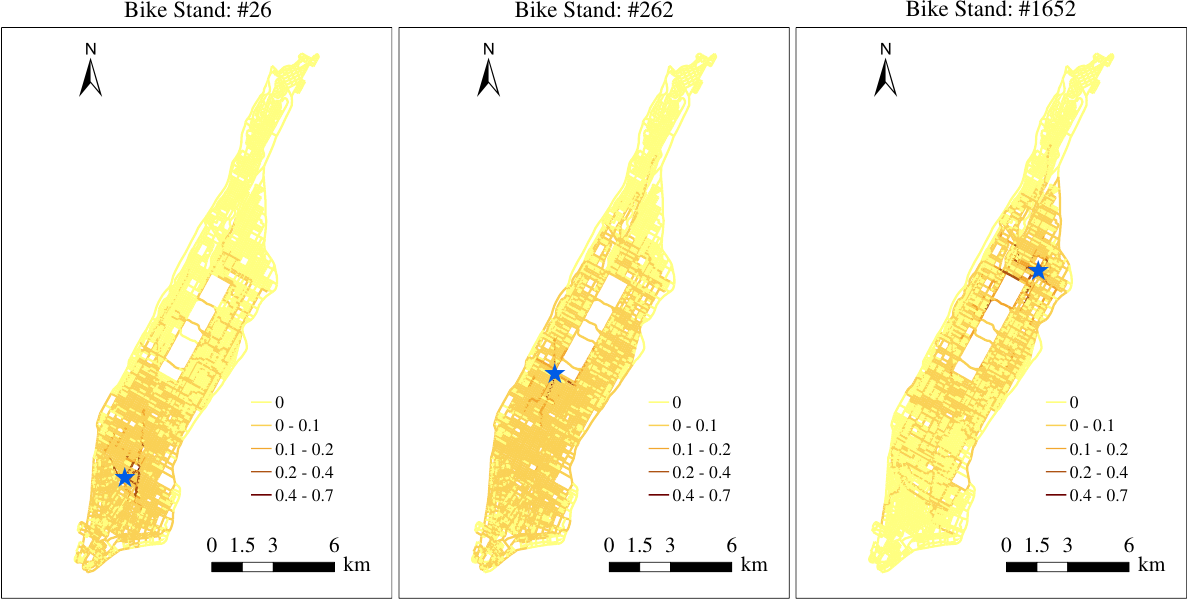}
		\caption{The visit probability $p_{s,e}$ of each bike from a fixed bike stand to different road sections. The blue five-pointed star represents the location of the bike stand, and the heatmap shows the values of the parameter $p_{s,e}$.}
		\label{fig_p_value}
	\end{figure}

	 Figure \ref{fig_p_value} visualizes the visit probability $p_{s,e}$ for bikes from a fixed bike stand to various road sections. The probability is higher for road segments closer to the bike stand, resulting in larger values of the fitting parameter $p_{s,e}$. This pattern reflects the operational characteristics of bike-sharing systems, where bikes tend to circulate more frequently in the vicinity of the bike stand, leading to increased coverage of nearby road segments.
	 
	\subsection*{Model of sensor-stand allocation}
	
	Given the total number of sensors $N_s$, this model allocates them across bike stands in a way that maximizes the number of road segments where the expected coverage meets or exceeds a predefined threshold, denoted as $K$. Let $n_s$ be the number of sensors allocated to bike stand $s$. Binary auxiliary variable $y_e=1$ if the expected coverage of road segment $e$ exceeds $K$ times. The MILP model of the sensor allocation problem is as follows:
	
	\begin{eqnarray}
		\label{SA_obj}
		\max_{n_s} \sum_{e \in \mathcal{A}}l_ey_e
	\end{eqnarray}
	
	\begin{eqnarray}
		\label{SA_1}
		 N_{e} = \sum_{s\in \mathcal{S}}p_{s,e}n_s  & &  \forall e \in \mathcal{A}
		\\
		\label{SA_2}
		-M(1-y_e) \leq  N_{e} - K  & & \forall e \in \mathcal{A}
		\\
		\label{SA_3}
		 N_{e} - K \leq My_{e}  & & \forall e \in \mathcal{A}
		\\
		\label{SA_4}
		\sum_{s \in \mathcal{S}}n_{s} \leq N_s
		\\
		\label{SA_5}
		0 \leq n_s \leq b_s  & & \forall s \in \mathcal{S} 
		\\
		n_s \in \mathbb{N}_{+}  & & \forall s \in \mathcal{S} 
		\\
		\label{SA_7}
		y_{e} \in \{0,1\}  & & \forall e \in \mathcal{A}
	\end{eqnarray}
	
	The objective function \eqref{SA_obj} maximizes the total length of road segments where the expected coverage meets the requirement. The expected coverage of road segment $e$ is determined by constraints \eqref{SA_1}. Constraints \eqref{SA_2} - \eqref{SA_3} express the binary state of coverage for road segment $e$, and are equivalent to the following via the big-M method:
	
    \begin{equation}\label{eqnconn}
	y_{e}=
	\begin{cases}
		1\quad &\text{if} ~~ N_{e} \geq K
		\\
		0 \quad &\text{if} ~~ N_{e} < K
	\end{cases}
	\end{equation}
	
	\noindent
	where $M$ is a sufficiently large number. Constraint \eqref{SA_4} ensures that the number of allocated sensors does not exceed the total. Constraints \eqref{SA_5} - \eqref{SA_7} are the decision variable constraints.
	
	In this study, we set $K=1$. The above model can be directly solved by off-the-shelf solvers.
	
	\subsection{Active scheduling of sensor-equipped bikes} \label{sec_active}
	In this section, we utilize the results obtained from the methods described in Section \ref{sec_MILP} to propose an active scheduling strategy specifically designed for sensor-equipped bikes to enhance sensing power further. The core idea is that when a user selects a bike at a stand, the app guides them to prioritize using a sensor-equipped bike. To capture real-world behavior, we introduce a guidance acceptance probability $\beta$, which simulates the randomness of user decisions by reflecting the likelihood that a user will follow the app's recommendation to choose a sensor-equipped bike. The detailed simulation process is outlined in Algorithm \ref{alg_active}.
 
	\begin{algorithm}[h!]
		\caption{Active scheduling of sensor-equipped bikes}
		\begin{tabbing}
			\hspace{0.01 in}\=  \hspace{0.9 in}\= \kill % set up two tab positions
			\>{\bf Input}  \> Set of bike trips $\mathcal{R}$, set of bike stands $\mathcal{S}$, simulation start time $t_0$, simulation 
			\\
			\>  \> end time $T$, initial number of bikes at each stand $b_{s}$, set of bikes at each 
			\\
			\> \> stand $\mathcal{K}_s, \forall s \in \mathcal{S}$,  number of sensors allocated to each bike stand $s$, probability
			\\
			\> \> $\beta$ of users accepting guidance strategies;
			\\
			\> {\bf Initialize} \> $t = t_0$, the status and type of each bike $k \in \mathcal{K}_s, \forall s \in \mathcal{S}$ is set to idle;
			\\
			\> {\bf 1:} \>  While $t \leq T$:
			\\
			\> {\bf 2:}\> \hspace{0.2 in} Check if any busy bikes need to end the trip at time $t$. If so, update their
			\\
			\> \>  \hspace{0.2 in} status to idle;
			\\
			\> {\bf 3:}\>  \hspace{0.2 in} Find all trips starting at time $t$, denoted as $\mathcal{R}^{'}$;
			\\
			\> {\bf 4:}\>   \hspace{0.2 in} For $r \in \mathcal{R}^{'}$:
			\\
			\> {\bf 5:}\>   \hspace{0.4 in} If $\text{Random}(0, 1) \leq \beta$ and there are idle sensor-equipped bikes at the
			\\
			\> \>   \hspace{0.4 in} start bike stand of trip $r$:
			\\
			\> {\bf 6:} \> \hspace{0.6 in} Randomly select an idle sensor-equipped bike at this stand;
			\\
			\> {\bf 7:}\>   \hspace{0.4 in} Else:
			\\
			\> {\bf 8:}\>   \hspace{0.6 in} Randomly select an idle bike from all bikes at this stand;
			\\
			\> {\bf 9:}\>   \hspace{0.4 in} End If
			\\
			\> {\bf 10:}\>  \hspace{0.4 in} Record the trajectory of the selected bike and update its status to busy;
			\\
			\> {\bf 11:}\> \hspace{0.2 in} End For
			\\
			\> {\bf 12:}\> \hspace{0.2 in} Update time: $t = t + 1$;
			\\
			\> {\bf 13:} \> End While
			\\
			\>{\bf Output}      \> The trajectory of each bike $k \in \mathcal{K}_s, \forall s \in \mathcal{S}$.
		\end{tabbing}
		\label{alg_active}
	\end{algorithm}
	
	\subsection{Quantify the sensing power of bike-sharing system} \label{sec_quantify}
	
	In the previous sections, we simulated the trajectories of each bike. In this section, we introduce a metric to quantify the sensing power of bike-sharing systems under different monitoring requirements.
	
	We divide the simulation time period into intervals $t \in \mathcal{T}$. The sensing reward is defined by the number of visits each road segment receives. Let $N_{e,t}$ represent the number of sensor-equipped bikes covering road segment $e$ during time interval $t$. The corresponding sensing reward $\xi(\cdot)$ is given by:
	
    \begin{equation}\label{eq_Q1}
	\xi(N_{e,t}) = 
	\begin{cases}
		1\quad &\text{if} ~~ N_{e,t} \geq 1
		\\
		0 \quad &\text{if} ~~ N_{e,t} < 1
	\end{cases}
	\end{equation}
	
	This binary sampling method tracks whether each spatial-temporal road segment meets the monitoring requirements of the application. By aggregating $N_{e,t}$ over all road segments and monitoring intervals, the overall sensing reward of the system is quantified as:
	
	\begin{equation}\label{eq_Q2}
	\Phi = \frac{\sum_{e \in \mathcal{A}}\sum_{t \in \mathcal{T}}l_{e}\xi(N_{e,t})}{N_{\mathcal{T}}\sum_{e \in \mathcal{A}}l_{e}} \times 100\%
	\end{equation}
	
	where $N_{\mathcal{T}}$ is the total number of time intervals. The numerator represents the total length of road segments meeting the sampling requirements, while the denominator reflects the potential total length if all segments meet the requirements. Therefore, $\Phi$ represents the proportion of spatial-temporal road segments covered.
	
	\section{Results and discussion} \label{sec_result}
	
	All the computational performances reported below are based on a Microsoft Windows 10 platform with Intel Core i9 - 3.60GHz and 16 GB RAM, using Python 3.8 and Gurobi 9.1.2.
	
	\subsection{Evaluation of the optimization strategies}
	In this section, we evaluate the effectiveness of sensor deployment strategies (Sec. \ref{sec_MILP}) and active scheduling (Sec. \ref{sec_active}) for enhancing the sensing power of the bike-sharing system. The three methods compared are as follows:
	
	\begin{itemize}
		\item \textbf{Random sensor allocation without active scheduling (Random-NoActive):} Sensors are randomly allocated to bike stands, and no active scheduling is applied to the sensor-equipped bikes.
		\item \textbf{Optimized sensor allocation without active scheduling (Optimized-NoActive):} Sensors are optimally allocated to stands using the method described in Sec. \ref{sec_MILP}, which maximizes spatial coverage. However, no active scheduling is applied.
		\item \textbf{Optimized sensor allocation with active scheduling (Optimized-Active):} In addition to optimized sensor allocation, sensor-equipped bikes are actively scheduled to enhance sensing performance using the method outlined in Sec. \ref{sec_active}, with $\beta = 1$ (assuming full user compliance with guidance).
	\end{itemize}
	
\begin{figure}[h!]
		\centering
		\includegraphics[width=0.9\textwidth]{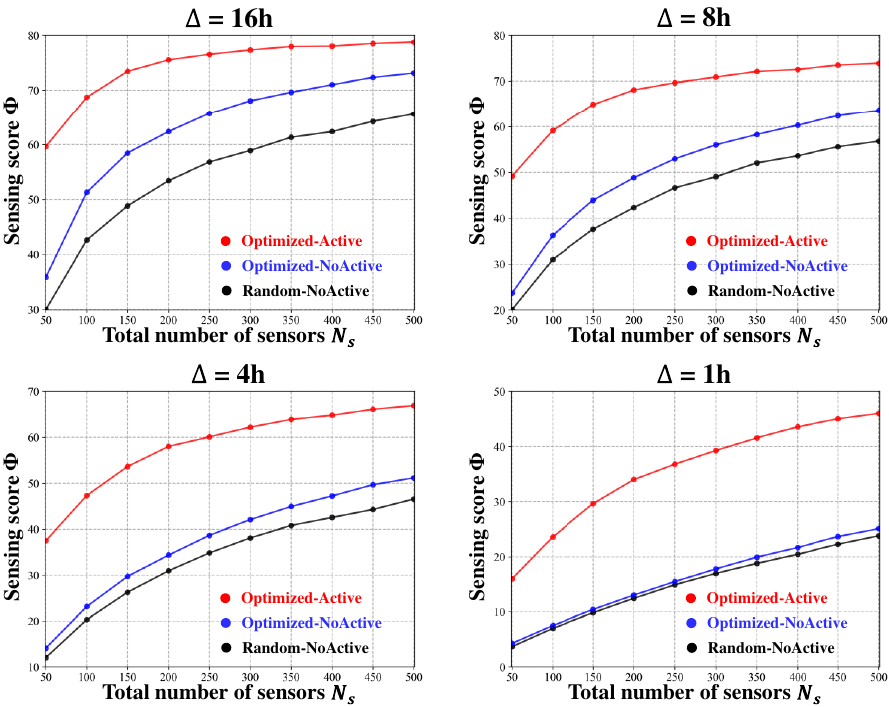}
		\caption{Sensing reward $\Phi$ for three methods across various total number of sensors $N_s$.}
		\label{fig_different_menthod}
	\end{figure}

	We plot the sensing score $\Phi$ generated by the above three methods for a given number of sensors ranging from 50 to 500. The results, shown in Figure \ref{fig_different_menthod}, are presented for four different sensing intervals: $\Delta = 16, 8, 4, 1$ (hours). The findings indicate that:
	\begin{itemize}
		\item[(1)] When the monitoring frequency demand is low (e.g. $\Delta = 16$), the optimized sensor allocation strategy increases the sensing reward by about 13\%-28\% compared to random sensor allocation. However, as the monitoring frequency demand increases (i.e., as $\Delta$ decreases), this advantage gradually diminishes.
		\item[(2)] The effectiveness of the optimized sensor allocation strategy becomes more evident with an increasing number of sensors, as it efficiently utilizes the available sensor resources.
		\item[(3)] The active scheduling strategy further enhances sensing power, particularly under high monitoring frequency requirements. For instance, when $\Delta = 4$ or $1$, active scheduling improves sensing reward by about 11.7\% - 24.1\% compared to scenarios without active scheduling, due to the improved utilization of sensor-equipped bikes through guided user behavior.
	\end{itemize}

\begin{figure}[h!]
		\centering
		\includegraphics[width=0.9\textwidth]{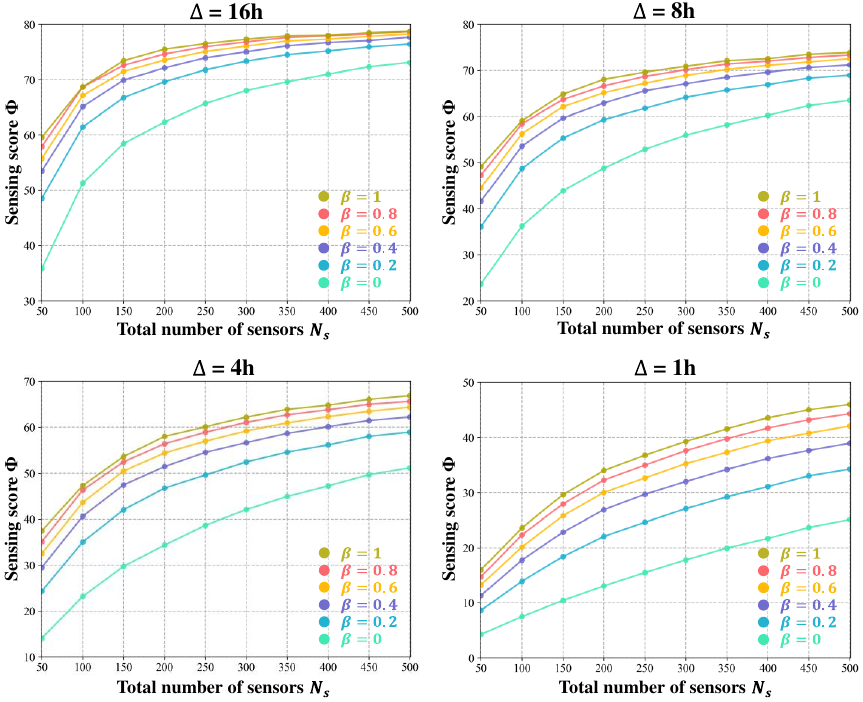}
		\caption{Sensing reward $\Phi$ for different guidance acceptance probabilities $\beta$ across various total number of sensors $N_s$.}
		\label{fig_different_beta}
	\end{figure}

	We further analyzed the impact of user guidance acceptance probability on sensing power, as shown in Figure \ref{fig_different_beta}. The results indicate that as $\beta$ increases from 0 to 0.2, the sensing reward $\Phi$ shows a significant rise. However, beyond $\beta = 0.6$, the marginal gain diminishes. This suggests that in practical applications, achieving guidance acceptance from 60\% of the users is sufficient for this strategy to be highly effective.

	\subsection{Sensing power of bike-sharing system}
	In this section, we evaluate the sensing power of the bike-sharing system. The trip data used in the subsequent analysis is sourced from Manhattan, covering the period from 6 am to 10 pm on March 1, 2024. This data has been processed using the methods outlined in Section \ref{sec_dp}. To serve these trips requires 6259 bikes by Algorithm \ref{alg_bike_num}. Additionally, the trajectories of sensor-equipped bikes were generated based on the optimized sensor allocation strategy detailed in Section \ref{sec_MILP} and the active scheduling strategy detailed in Section \ref{sec_active}. 
	
%	\begin{figure}[h!]
%		\centering
%		\includegraphics[width=0.8\textwidth]{fig_trip_num.pdf}
%		\caption{Number of bike-sharing trips starting each hour.}
%		\label{fig_trip_num}
%	\end{figure}
	
\begin{figure}[h!]
	 	\centering
	 	\includegraphics[width=0.8\textwidth]{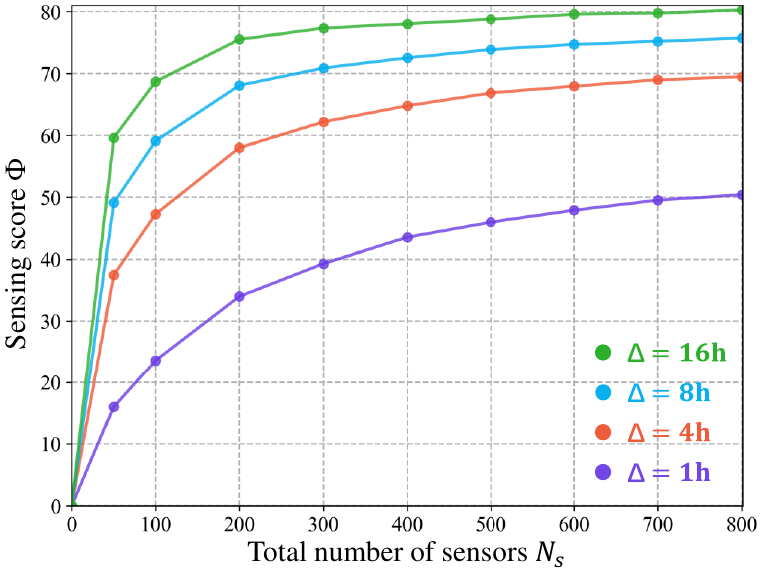}
	 	\caption{Sensing reward $\Phi$ for different sensing intervals $\Delta$ across various total number of sensors $N_s$.}
	 	\label{fig_sensing_power}
 	\end{figure}

	Figure \ref{fig_sensing_power} presents the sensing reward $\Phi$ for different sensing intervals $\Delta$ across various total number of sensors $N_s$. The results indicate that: 
	\begin{itemize}
		\item[(1)] \textbf{Significant sensing capabilities for low-frequency monitoring:} The bike-sharing system exhibits substantial sensing capabilities, which makes it highly suitable for low-frequency monitoring scenarios such as built environment and urban greenery health monitoring. For example, equipping sensors on just 100 bikes (approximately 1\% of the total bike fleet in Manhattan) achieves coverage of nearly 70\% of road segments within a single day (16 hours). This shows that even with a limited number of sensors, extensive coverage can be achieved for applications where frequent data collection is not critical.
		
		\item[(2)] \textbf{Sensor requirements for different monitoring frequencies:} Achieving 50\% spatiotemporal coverage requires varying numbers of sensors depending on the monitoring frequency. High-frequency applications, such as air quality monitoring and urban heat island studies, which demand detailed and frequent data, require about 800 sensors for a 1-hour monitoring interval ($\Delta=1$). In contrast, applications that can operate with less frequent monitoring, such as those with 4-hour intervals ($\Delta=4$), need only 121 sensors, while 8-hour ($\Delta=8$) and 16-hour ($\Delta=16$) intervals require 54 and 41 sensors, respectively. These results highlight a sharp increase in sensor requirements as the monitoring frequency intensifies. Furthermore, while high-frequency applications necessitate a larger number of sensors, complete data capture may not always be essential, and data gaps can be addressed using appropriate inference models.
		
		\item[(3)] \textbf{Diminishing returns with increased sensors:}  As the number of sensors increases, the marginal gain in sensing reward decreases, indicating diminishing returns. This suggests that beyond a certain point, adding more sensors provides limited additional benefits.
	\end{itemize}

\begin{figure}[h!]
		\centering
		\includegraphics[width=0.8\textwidth]{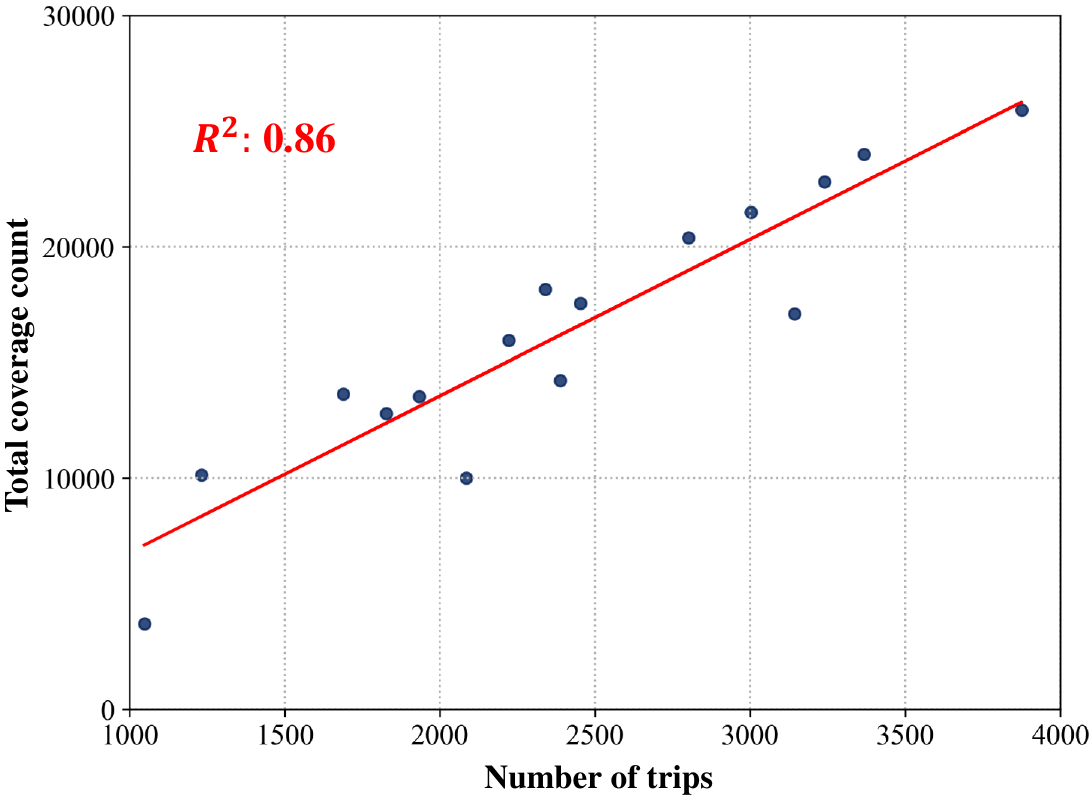}
		\caption{Scatter plot demonstrating the relationship between the hourly number of trips and the total coverage of road segments by sensor-equipped bikes.}
		\label{fig_cover_count}
	\end{figure}

	In addition to the overall sensing reward, we further analyze the bike-sharing system's sensing power from an hourly spatial perspective.
	Taking the deployment of 800 sensor-equipped bikes as an example, we first calculated the total coverage count of these sensor-equipped bikes on each road segment for every hour. Figure \ref{fig_cover_count} presents a scatter plot, with the x-axis representing the number of trips in that hour and the y-axis representing the total coverage count across all road segments. The results indicate a positive correlation between the number of trips and the total road segment coverage, suggesting that increased bike activity enhances the system's sensing capabilities.
	
\begin{figure}[h!]
		\centering
		\includegraphics[width=1.0\textwidth]{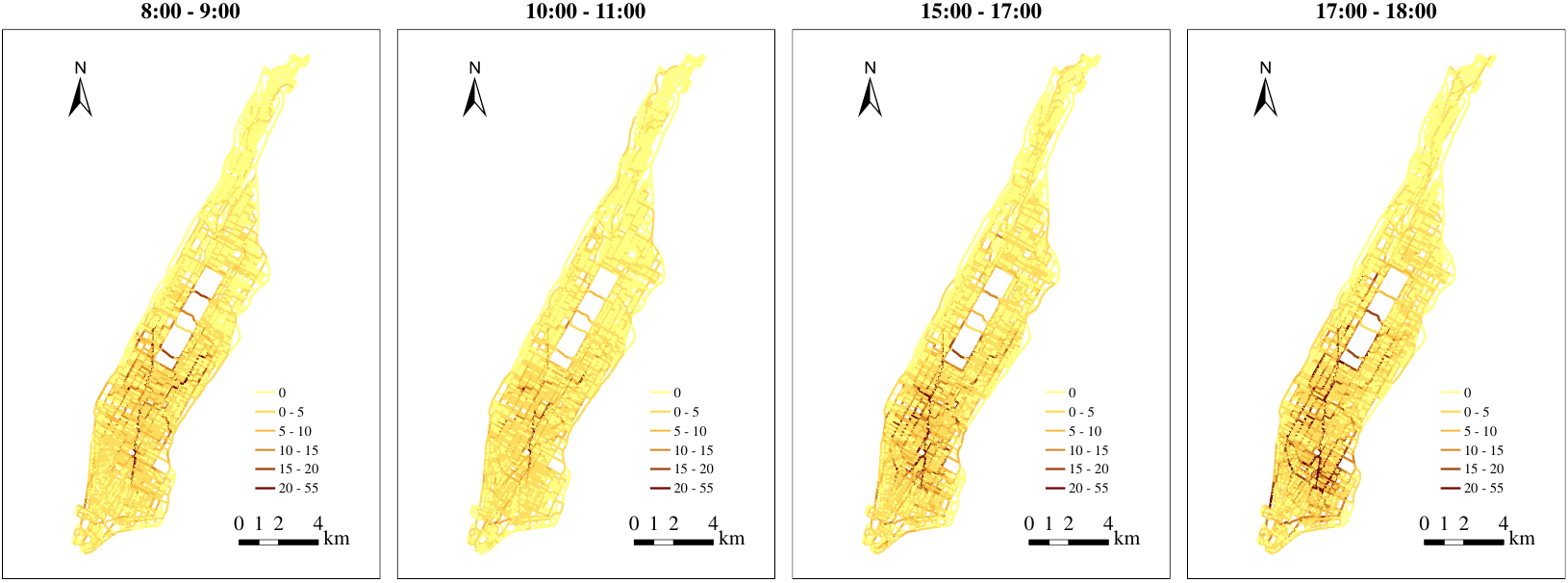}
		\caption{Heatmap of road segment coverage by 800 sensor-equipped bikes over various hours.}
		\label{fig_spatial}
	\end{figure}

	Next, we created heatmaps to illustrate the coverage counts of road segments by sensor-equipped bikes during various hours, as shown in Figure \ref{fig_spatial}. Our analysis reveals a consistent spatial distribution of road network coverage across different time intervals. Notably, the southern region exhibits a higher frequency of scanned road segments, while the northern region shows significantly less coverage. This pattern suggests that certain road segments consistently receive more attention, likely correlating with areas of increased bike activity and trip density. Furthermore, the heatmaps indicate that the spatial configuration of bike usage remains relatively stable throughout the day. These insights can guide strategies for optimizing sensor deployment and scheduling, ensuring that high-traffic areas receive adequate coverage while addressing potential gaps in less frequented regions. 
	
	To address the road segments that are rarely scanned, we can implement fixed monitoring stations or deploy dedicated sensing fleets. These strategies would ensure that underrepresented areas receive the necessary attention for data collection, enhancing the overall sensing capability of the system. By strategically placing fixed monitoring stations in key locations and utilizing dedicated sensing vehicles to target these less frequented segments, we can fill in coverage gaps and improve the comprehensiveness of the spatial data gathered by the bike-sharing system. 
	
	\section{Conclusion} \label{sec_conclusion}
	This study leverages bike-sharing systems as a cost-effective urban sensing platform by introducing strategies to optimize sensor deployment and actively schedule sensor-equipped bikes to monitor complex urban dynamics. A case study in Manhattan demonstrates the effectiveness of these strategies and quantifies the sensing power of bike-sharing systems, highlighting their significant potential as mobile sensing networks. The key findings from this study are as follows:
	
	\begin{itemize}
		\item The optimized sensor allocation strategy enhances sensing rewards by approximately 13\%-28\% compared to random sensor allocation when monitoring once per road segment per day.
		\item In practical applications, once more than 60\% of users follow the guided bike selection, there is a significant improvement in sensing rewards. However, further increases beyond this 60\% acceptance rate yield only marginal gains.
		\item Deploying sensors on approximately 1\% of the bikes in Manhattan achieves coverage of about 70\% of road segments within a 16-hour period.
		\item  The required number of sensors escalates substantially with higher monitoring frequency demands.
		\item As the number of sensors increases, the marginal gain in coverage decreases, indicating that sensor deployment strategies should carefully consider the trade-off between coverage and resource allocation.
		\item There is a strong consistency in the spatial coverage of sensor-equipped vehicles across different hours, suggesting that road segments frequently scanned at one time continue to receive attention throughout the day. This insight facilitates the strategic placement of fixed monitoring stations or the deployment of dedicated vehicles in areas that are challenging to cover.
	\end{itemize}
	
	Our study fills a critical gap in current transportation research, demonstrating that bike-sharing systems offer a sustainable, flexible, and efficient solution for urban sensing. The findings have significant implications for policymakers, urban planners, and transportation professionals seeking to implement data-driven strategies for monitoring environmental and urban dynamics. By leveraging bike-sharing networks, cities can gain extensive spatial-temporal insights into urban environments, ultimately contributing to smarter, more resilient urban planning and management.
	
	Future explorations can focus on the following areas: First, developing real-time, adaptive scheduling and dynamic deployment strategies for sensor-equipped bikes is essential to address the changing monitoring demands across different times and regions. Additionally, integrating bike-sharing systems with other types of mobile sensors, such as those mounted on taxis or buses, can enhance overall spatial-temporal sensing power, paving the way for a more comprehensive urban sensing network.
	
	\section*{Declaration of competing interest}
	The authors declare no conflict of interests.
	
	\section*{Acknowledgements}
	This work was supported by the National Natural Science Foundation of China [grant numbers 62394331, 62394335, 52178044, 72071163, 72101215].


\begin{thebibliography}{99}
	\bibitem[Agarwal et al., 2020]{Agarwal2020} Agarwal, D. , Iyengar, S. , Swaminathan, M. , Sharma, E. , Raj, A. , Hatwar, A., 2020. Modulo: Drive-by Sensing at City-scale on the Cheap, in: Proceedings of the 3rd ACM SIGCAS Conference on Computing and Sustainable Societies. Ecuador, 187-197. 
	
	\bibitem[Ali et al., 2021]{AQSS2021} Ali, A., Qureshi, M.A., Shiraz, M., Shamim, A., 2021. Mobile Crowd Sensing Based Dynamic Traffic Efficiency Framework for Urban Traffic Congestion Control. Sustainable Computing: Informatics and Systems 32, 100608. 
	
	\bibitem[Alsina-Pag\'es et al., 2017]{Alsina-Pages2017} Alsina-Pag\'es, R., Hernandez-Jayo, U., Alías, F., Angulo, I., 2017. Design of a Mobile Low-Cost Sensor Network Using Urban Buses for Real-Time Ubiquitous Noise Monitoring. Sensors 17 (1), 57.
	
	\bibitem[Anjomshoaa et al., 2018]{ADRMdR2018} Anjomshoaa, A., Duarte, F., Rennings, D., Matarazzo, T.J., deSouza, P., Ratti, C., 2018. City Scanner: Building and Scheduling a Mobile Sensing Platform for Smart City Services. IEEE Internet of Things Journal, 5 (6), 4567-4579.
	
	\bibitem[Ariss et al., 2024]{AWSDR2024} Ariss, M., Wang, A., Sabouri, S., Duarte, F., Ratti, C., 2024. Drive-by Environmental Sensing Strategy to Reach Optimal and Continuous Spatio-Temporal Coverage Using Local Transit Network. Transportation Research Record.
	
	\bibitem[Asprone et al., 2021]{ADFS2021} Asprone, D., Di Martino, S., Festa, P., Starace, L.L.L., 2021. Vehicular Crowd-sensing: a Parametric Routing Algorithm to Increase Spatio-temporal Road Network Coverage.  International Journal of Geographical Information Science, 35 (9), 1876-1904.
	
	\bibitem[Birenboim et al., 2021]{BHK2021} Birenboim, A., Helbich, M., Kwan, M.-P., 2021. Advances in Portable Sensing for Urban Environments: Understanding Cities from a Mobility Perspective. Computers, Environment and Urban Systems, 88, 101650.
	
	\bibitem[Chen et al., 2024]{CMLZH2024} Chen, Q., Ma, S., Li, H., Zhu, N., He, Q., 2024. Optimizing Bike Rebalancing Strategies in Free-floating Bike-sharing Systems: An Enhanced Distributionally Robust Approach. Transportation Research Part E: Logistics and Transportation Review, 184, 103477. 
	
	\bibitem[Chen et al., 2020]{C2020} Chen, X., Xu, S., Han, J., Fu, H., Pi, X., Joe-Wong, C., Li, Y., Zhang, L., Noh, H.Y., Zhang, P., 2020. PAS: Prediction-Based Actuation System for City-Scale Ridesharing Vehicular Mobile Crowdsensing. IEEE Internet of Things Journal, 7 (5), 3719-3734.
	
	\bibitem[Dai and Han, 2023]{DH2023} Dai, Z., Han, K., 2023. Exploring the Drive-by Sensing Power of Bus Fleet through Active Scheduling. Transportation Research Part E: Logistics and Transportation Review, 171, 103029.
	
	\bibitem[Du et al., 2019]{DSXVF2019} Du, R., Santi, P., Xiao, M., Vasilakos, A.V., Fischione, C., 2019. The Sensable City: A Survey on the Deployment and Management for Smart City Monitoring. IEEE Communications Surveys \& Tutorials 21 (2), 1533-1560.
	
	\bibitem[Fekih et al., 2021]{FBRDRAP2021} Fekih, M.A., Bechkit, W., Rivano, H., Dahan, M., Renard, F., Alonso, L., Pineau, F., 2021. Participatory Air Quality and Urban Heat Islands Monitoring System. IEEE Transactions on Instrumentation and Measurement 70, 1-14.
	
	\bibitem[Gao et al., 2016]{Gao2016} Gao, Y., Dong, W., Guo, K., Liu, Xue, Chen, Y., Liu, Xiaojin, Bu, J., Chen, C., 2016. Mosaic: A Low-cost Mobile Sensing System for Urban Air Quality Monitoring. in: IEEE INFOCOM 2016 - The 35th Annual IEEE International Conference on Computer Communications, San Francisco, CA, USA, 1-9. 
	
	\bibitem[Glock and Meyer, 2020]{GM2020} Glock, K., Meyer, A., 2020. Mission Planning for Emergency Rapid Mapping with Drones. Transportation Science, 54, 534–560. 
	
	\bibitem[Guo and Qian, 2024]{GQ2024} Guo, S., Qian, X., 2024. Optimal Drive-By Sensing in Urban Road Networks With Large-Scale Ridesourcing Vehicles. IEEE Transactions on Intelligent Transportation Systems. 1–12.
	
	\bibitem[Gupta et al., 2024]{GMPDPR2024} Gupta, A., Mora, S., Preisler, Y., Duarte, F., Prasad, V., Ratti, C., 2024. Tools and Methods for Monitoring the Health of the Urban Greenery. Nature Sustainability. 7, 536–544.
	
	
	
	\bibitem[Han et al., 2024]{HJNLL2024} Han, K., Ji, W., Nie, Y., Li, Z., Liu, S., 2024. Exploring the Sensing Power of Mixed Vehicle Fleets. Transportation Research Part B: Methodological, 190, 103066.
	
	\bibitem[Hasenfratz et al., 2015]{HSWHFABT2015} Hasenfratz, D., Saukh, O., Walser, C., Hueglin, C., Fierz, M., Arn, T., Beutel, J., Thiele, L., 2015. Deriving High-resolution Urban Air Pollution Maps Using Mobile Sensor Nodes. Pervasive and Mobile Computing 16, 268-285.
	
	\bibitem[Huang et al., 2024]{HLYK2024} Huang, M., Li, X., Yang, M., Kuai, X., 2024. Intelligent Coverage and Cost-effective Monitoring: Bus-based Mobile Sensing for City Air Quality. Computers, Environment and Urban Systems, 108, 102073.
	
	\bibitem[Hua et al., 2024]{HYCCC2024} Hua, M., Yu, X., Chen, X., Chen, J., Cheng, L., 2024. Can Bike Sharing Achieve Self-balancing Distribution? Evidence from Dockless and Station-based Cases. Travel Behaviour and Society, 38, 100879.
	
	
	
	\bibitem[Ji et al., 2023a]{JHL2023a} Ji, W., Han, K., Liu, T., 2023a. A Survey of Urban Drive-by Sensing: An Optimization Perspective. Sustainable Cities and Society, 99, 104874.
	
	\bibitem[Ji et al., 2023b]{JHL2023b} Ji, W., Han, K., Liu, T., 2023b. Trip-based Mobile Sensor Deployment for Drive-by Sensing with Bus Fleets. Transportation Research Part C: Emerging Technologies, 157, 104404.
	
	\bibitem[Ji et al., 2023c]{JHG2023} Ji, W., Han, K., Ge, Q., 2023c. Route Planning of Mobile Sensing Fleets for Repeatable Visits, arXiv: 2307.02397. 
	
	\bibitem[Kaivonen and Ngai, 2020]{Kaivonen2020} Kaivonen, S., Ngai, E.C.-H., 2020. Real-time Air Pollution Monitoring with Sensors on City Bus. Digital Communications and Networks, 6 (1), 23-30.
	
	\bibitem[Kon et al., 2022]{KFDDSR2022} Kon, F., Ferreira, É.C., De Souza, H.A., Duarte, F., Santi, P., Ratti, C., 2022. Abstracting Mobility Flows from Bike-sharing Systems. Public Transport, 14, 545–581.
	
	\bibitem[Kong et al., 2020]{KJS2020} Kong, H., Jin, S.T., Sui, D.Z., 2020. Deciphering the Relationship Between Bikesharing and Public Transit: Modal Substitution, Integration, and Complementation. Transportation Research Part D: Transport and Environment, 85, 102392.
	
	\bibitem[Li et al., 2018]{LZB2018} Li, P., Zhao, P., Brand, C., 2018. Future Energy Use and CO2 Emissions of Urban Passenger Transport in China: A Travel Behavior and Urban form Based Approach. Applied Energy, 211, 820–842. 
	
	\bibitem[Lin et al., 2018]{LHP2018} Lin, L., He, Z., Peeta, S., 2018. Predicting Station-level Hourly Demand in a Large-scale Bike-sharing Network: A Graph Convolutional Neural Network Approach. Transportation Research Part C: Emerging Technologies, 97, 258–276. 
	
	\bibitem[Li and Long, 2024]{LL2024} Li, Y., Long, Y., 2024. Inferring Storefront Vacancy Using Mobile Sensing Images and Computer Vision Approaches. Computers, Environment and Urban Systems 108, 102071. 
	
	\bibitem[Li et al., 2023]{LMZLL2023} Li, Y., Meng, X., Zhao, H., Li, W., Long, Y., 2023. Identifying abandoned buildings in shrinking cities with mobile sensing images. Urban Informatics, 2(1), 3.

	\bibitem[Ma et al., 2014]{MZY2014} Ma, H., Zhao, D., Yuan, P., 2014. Opportunities in Mobile Crowd Sensing. IEEE Communications Magazine 52 (8), 29-35.
	
	\bibitem[Meng and Han, 2023]{MH2023} Meng, D., Han, K., 2023. Optimizing Vehicle-Passenger Matching for Online Ride-Hailing with Vehicular Crowd-Sensing, in: 2023 IEEE 26th International Conference on Intelligent Transportation Systems (ITSC). Presented at the 2023 IEEE 26th International Conference on Intelligent Transportation Systems (ITSC), IEEE, Bilbao, Spain, pp. 3527–3532.
	
	\bibitem[Messier et al., 2018]{Messier2018} Messier, K.P., Chambliss, S.E., Gani, S., Alvarez, R., Brauer, M., Choi, J.J., Hamburg, S.P., Kerckhoffs, J., LaFranchi, B., Lunden, M.M., Marshall, J.D., Portier, C.J., Roy, A., Szpiro, A.A., Vermeulen, R.C.H., Apte, J.S., 2018. Mapping Air Pollution with Google Street View Cars: Efficient Approaches with Mobile Monitoring and Land Use Regression. Environmental Science \& Technology, 52 (21), 12563-12572.
	
	\bibitem[O'Keeffe et al., 2019]{OASSR2019} O'Keeffe, K.P., Anjomshoaa, A., Strogatz, S.H., Santi, P., Ratti, C., 2019. Quantifying the Sensing Power of Vehicle Fleets. Proceedings of the National Academy of Sciences, 116 (26), 12752–12757.
	
	
	\bibitem[Shen et al., 2022]{SLZSW2022} Shen, S., Lv, C, Zhu, H., Sun, L., Wang, R., 2022. Potentials and Prospects of Bicycle Sharing System in Smart Cities: A Review. IEEE Sensors Journal, 22, 7519–7533.
	
	\bibitem[SM et al., 2019]{SM2019} SM, S. N., Yasa, P. R., Narayana, M. V., Khadirnaikar, S., Rani, P., 2019. Mobile Monitoring of Air Pollution Using Low Cost Sensors to Visualize Spatio-temporal Variation of Pollutants at Urban Hotspots. Sustainable Cities and Society, 44, 520-535. 
	
	\bibitem[Voinea et al., 2020]{VBP2020} Voinea, S.C., Bujari, A., Palazzi, C.E., 2020. Air Quality Control through Bike Sharing Fleets, in: 2020 IEEE Symposium on Computers and Communications (ISCC). Presented at the 2020 IEEE Symposium on Computers and Communications (ISCC), IEEE, Rennes, France, pp. 1–4.
	
	\bibitem[Wu et al., 2020]{WXLLWZLG2020} Wu, D., Xiao, T., Liao, X., Luo, J., Wu, C., Zhang, S., Li, Y., Guo, Y., 2020. When Sharing Economy Meets IoT: Towards Fine-grained Urban Air Quality Monitoring through Mobile Crowdsensing on Bike-share System. Proc. ACM Interact. Mob. Wearable Ubiquitous Technol. 4, 1–26. 
	
	\bibitem[Xu et al., 2018]{XJL2018} Xu, C., Ji, J., Liu, P., 2018. The Station-free Sharing Bike Demand Forecasting with a Deep Learning Approach and Large-scale Datasets. Transportation Research Part C: Emerging Technologies, 95, 47–60. 
	
	\bibitem[Xu et al., 2019]{XCPJZN} Xu, S., Chen, X., Pi, X., Joe-Wong, C., Zhang, P., Noh, H.Y., 2019. iLOCuS: Incentivizing Vehicle Mobility to Optimize Sensing Distribution in Crowd Sensing. IEEE Transactions on Mobile Computing. 19, 1831–1847. 
	
	
	
	\bibitem[Yang et al., 2019]{YHTC2019} Yang, Y., Heppenstall, A., Turner, A., Comber, A., 2019. A Spatiotemporal and Graph-based Analysis of Dockless Bike Sharing Patterns to Understand Urban Flows Over the Last Mile. Computers, Environment and Urban Systems, 77, 101361.
	
	\bibitem[Yu et al., 2022]{YXLZLSCYY2022} Yu, Q., Xie, Y., Li, W., Zhang, H., Liu, X., Shang, W., Chen, J., Yang, D., Yan, J., 2022. GPS Data in Urban Bicycle-sharing: Dynamic Electric Fence Planning with Assessment of Resource-saving and Potential Energy Consumption Increasement. Applied Energy, 322, 119533. 
	
	\bibitem[Yu et al., 2024]{YHLLCKL2024} Yu, S., Han, X., Liu, L., Liu, G., Cheng, M., Ke, Y., Li, L., 2024. Exploring usage pattern variation of free-floating bike-sharing from a night travel perspective. Scientific reports, 14(1), 16017.
	
	\bibitem[Zhang et al., 2019]{ZSLXFZHSL2019} Zhang, H., Song, X., Long, Y., Xia, T., Fang, K., Zheng, J., Huang, D., Shibasaki, R., Liang, Y., 2019. Mobile Phone GPS Data in Urban Bicycle-sharing: Layout Optimization and Emissions Reduction Analysis. Applied Energy, 242, 138–147. 
	
	\bibitem[Zhang et al., 2018]{ZYL2018} Zhang, X., Yang, Z., Liu, Y., 2018. Vehicle-based Bi-objective Crowdsourcing. IEEE Transactions on Intelligent Transportation Systems, 19 (10), 3420-3428.
	
	\bibitem[Zhang et al., 2019]{ZLL2019} Zhang, Y., Lin, D., Liu, X. C., 2019. Biking islands in cities: An analysis combining bike trajectory and percolation theory. Journal of Transport Geography, 80, 102497. 
	
	\bibitem[Zhao et al., 2014]{ZML2014} Zhao, D., Ma, H., Liu, L., 2014. Energy-efficient Opportunistic Coverage for People-centric Urban Sensing. Wireless Networks, 20 (6), 1461-1476. 
		
	\end{thebibliography}
\end{document}